\documentclass{amsart}

\usepackage{amsmath}
\usepackage{graphics}

\newtheorem{Thm}{Theorem}
\newtheorem{Prop}[Thm]{Proposition}

\newtheorem*{Thm*}{Theorem}

\newcommand{\lora}{{\longrightarrow}}

\newcommand{\de}{\partial}
\newcommand{\R}{\mathbb{R}}

\newcommand{\imbr}[1]{{\,\hbox{\rm Imb}\,(S^1,\R^{#1}) }}

\begin{document}

\markboth{Riccardo Longoni}
{Nontrivial classes in $H^*(\imbr n)$ from nontrivalent graph cocycles}

\title{Nontrivial classes in $H^*(\imbr n)$ from nontrivalent graph cocycles}

\author{Riccardo Longoni}

\address{Dipartimento di Matematica ``G. Castelnuovo'',
Universit\`a di Roma ``La Sapienza'',\\ Piazzale Aldo Moro, 2 --
I-00185 Roma, Italy}
\email{longoni@mat.uniroma1.it}

\maketitle

\begin{abstract}
We construct nontrivial cohomology classes of the space $\imbr n$ 
of imbeddings of the circle into $\R^n$, by means of Feynman diagrams. 
More precisely, starting from a suitable linear combination of nontrivalent 
diagrams, we construct, for every even number $n\geq 4$, a de~Rham cohomology 
class on $\imbr n$. We prove nontriviality of these classes by evaluation 
on the dual cycles.
\end{abstract}

\section{Introduction}

In the recent years Quantum Field Theories have provided useful tools 
for addressing problems in Algebraic Topology. An example is
the computation of the real cohomology of the spaces $\imbr n$ of 
imbeddings of the circle into $\R^n$ for $n\ge 3$. The works of Bar-Natan 
\cite{BN95} and Kontsevich \cite{K}, have shown how the perturbative 
expansion of Chern-Simons Quantum Field Theory yields the so-called 
``Vassiliev invariants'' of the space of ordinary knots. The precise 
statement is that there exists a complex of (trivalent) Feynman diagrams 
whose cohomology classes give rise to 0-cohomology classes of the space 
$\imbr 3$.

Chern-Simons is a Topological Quantum Field Theory 
which can be defined in three dimensions only and produces Feynman 
diagrams whose vertices have valence equal to three. There exist however 
other Quantum Field Theories that, while sharing the same topological 
properties, are defined over manifolds of any dimension and produce
Feynman diagrams not necessarily trivalent. These Topological Quantum Field 
Theories are known as ``BF Theories'' and their properties are studied in 
details in \cite{CCRRo} and \cite{CCRo}.

In \cite{CCRL} it is shown how the perturbative expansion of 
BF Theories can define a complex of Feynman diagrams and a chain map into 
the de~Rham complex of $\imbr n$, for any $n>3$. 
Moreover this map gives rise to an 
injective map in cohomology, when restricted to trivalent diagrams. 
The cohomology of this complex of Feynman diagrams is called 
``graph cohomology''.

A different approach to the problem of computing the cohomology of 
$\imbr n$ has been considered by V. Vassiliev in \cite{V}: there exists in 
fact a spectral sequence converging to the cohomology of $\imbr n$ for every 
$n\geq 4$, whose $E_1$ term coincides with the graph cohomology (see also
\cite{T1}).
Thus, Vassiliev spectral sequence gives an ``upper bound'' to the cohomology 
of $\imbr n$, since nontrivial cohomology classes
of $\imbr n$ must come from some element in $E_1$, while \cite{CCRL} gives 
a ``lower bound'', since at least those elements corresponding to trivalent 
diagrams give rise to nontrivial elements in cohomology. In general it is 
unclear whether the collapsing always happens at the $E_1$ term, 
as conjectured 
by Kontsevich \cite{K}, or there are some ``graph cocycles'' which do not 
give rise to classes in $H^*(\imbr n)$. 

In this note we want to extend the results of \cite{CCRL} in the following 
direction
\begin{Thm}
\label{thm-main}
For every even $n\geq 4$, there exists a nontrivial class in\break
$H^{(n-3)3+1}(\imbr n)$ associated to a nontrivalent graph cocycle.
\end{Thm}

A similar question has been posed by R. Bott \cite{Bott} regarding the
possibility of constructing 1-dimensional cocycles on the space $\imbr
3$ of ordinary knots. In this case, however, it is not even established
the convergence of the spectral sequence. We also mention the work of
D. Sinha \cite{Si} on the cohomology of the spaces of imbeddings of the 
real line into $\R^n$. Using Goodwillie's calculus, one can define a spectral 
sequence whose $E_1$-term, as in Vassiliev's case, is isomorphic to (a 
variant of) the graph cohomology. Recently P. Lambrechts and I. Volic
have announced a proof of the collapse of Sinha spectral sequence at the 
$E_1$ term.\\

Here is the plan of the paper. 
We will first recall the definition of graph cohomology
(Sections~\ref{sec:grco} and \ref{sec:grho}), and the construction of
the map from graph cohomology to the de~Rham cohomology of $\imbr n$
(Section~\ref{sec:csi}). Then, in Section~\ref{sec:nontriv}, we will
prove Theorem~\ref{thm-main} by considering a nontrivial graph cocycle
consisting of nontrivalent diagrams and mapping it to $H^*(\imbr n)$. 
Nontriviality of this class will follow from an explicit evaluation on 
the dual cycle.


\section{Graph cohomology}
\label{sec:grco}

We first briefly recall some definitions given in \cite{CCRL}. By
$\mathcal D_o^{k,m}$ and $\mathcal D_e^{k,m}$ we mean the real vector
spaces generated by decorated diagrams of order $k$ and degree $m$, of
odd and even type, respectively. The diagrams consist of an {\em
oriented circle}\/ and many {\em edges}\/ joining vertices which may
lie either on the circle ({\em external vertices}) or off the circle
({\em internal vertices}). We require all the vertices to be at 
least trivalent. If we denote by $e$ the number of edges of a diagram, 
$v_i$ the number of its
internal vertices and $v_e$ the number of its external vertices, then
the order $k$ of a diagram is $e-v_i$ (minus the Euler
characteristic) while the degree $m$ is $2e-3v_i-v_e$ (the
deviation from being a trivalent diagram).  
We also define a {\em chord diagram} to be a decorated diagram whose
vertices are all external. 

The decoration in $\mathcal D_o^{k,m}$ is given by numbering the
vertices (up to even cyclic permutations), numbering the
internal vertices (up to even permutations) and orienting the edges
(up to reversal). By convention, we number external vertices from 1
to $v_i$ and internal vertices from $v_i+1$ to $v_i+v_e$. 
An extra decoration is needed on the edges connecting
the same external vertex, namely an ordering of the two half-edges
forming them. The decoration in $\mathcal D_e^{k,m}$ is given by
numbering the external vertices (up to even cyclic permutations) and
numbering the edges (up to even permutations). Finally, both in
$\mathcal D_o^{k,m}$ and $\mathcal D_e^{k,m}$ we quotient by the
subspace generated by all diagrams containing two edges joining the
same pair of vertices and diagrams containing edges whose end-points 
are the same internal vertex.

The coboundary operators $\delta_o\colon \mathcal
D^{k,m}_o \to \mathcal D^{k,m+1}_o$ and $\delta_e\colon \mathcal
D^{k,m}_e \to \mathcal D^{k,m+1}_e$ are linear operators whose action
on a diagram $\Gamma$ is given by the signed sum of all the diagrams obtained
from $\Gamma$ by contracting, one at a time, all the regular edges and 
arcs of the diagrams. Here by arc we mean a piece of the oriented circle 
between two consecutive vertices, and by regular edge we mean an edge with 
at least one internal end-point. The signs are as follows: if we contract
an arc or edge connecting the vertex $i$ with the vertex $j$, and oriented
from $i$ to $j$, then the sign is $(-1)^{j}$ if $j>i$ or $(-1)^{i+1}$
if $j<i$. If we contract an edge labelled by $\alpha$, then the
sign is $(-1)^{\alpha+1+v_e}$,where $v_e$ is the number of external
vertices of the diagram. 
 
These complexes are called {\em graph complexes}, and the
relevant cohomology groups $H^{k,m}(\mathcal{D}_o)$ and
$H^{k,m}(\mathcal{D}_e)$ are know as {\em graph cohomology}. When
we write $(\mathcal D^{k,m},\delta)$ (resp. $H^{k,m}(\mathcal{D})$), we
mean either the odd or the even-type graph complex (resp. graph
cohomology). 


\section{Graph homology}
\label{sec:grho}

Graph homology is given by the dual vector
space $(\mathcal D^{k,m})^*$ and the boundary operator $\de\colon
(\mathcal D^{k,m})^* \to (\mathcal D^{k,m-1})^*$ defined as the
adjoint to $\delta$. The homology groups are denoted by
$H_{k,m}(\mathcal D)$. Notice that since there is a
preferred basis of $\mathcal D^{k,m}$ given by the diagrams, we can
identify $(\mathcal D^{k,m})^*$ with $\mathcal D^{k,m}$, and
consider $\de$ as an operator on $\mathcal D^{k,m}$ that
decollapses the vertices of the diagrams in all possible ways.

\begin{figure}[h]
\resizebox{5in}{!}{\includegraphics{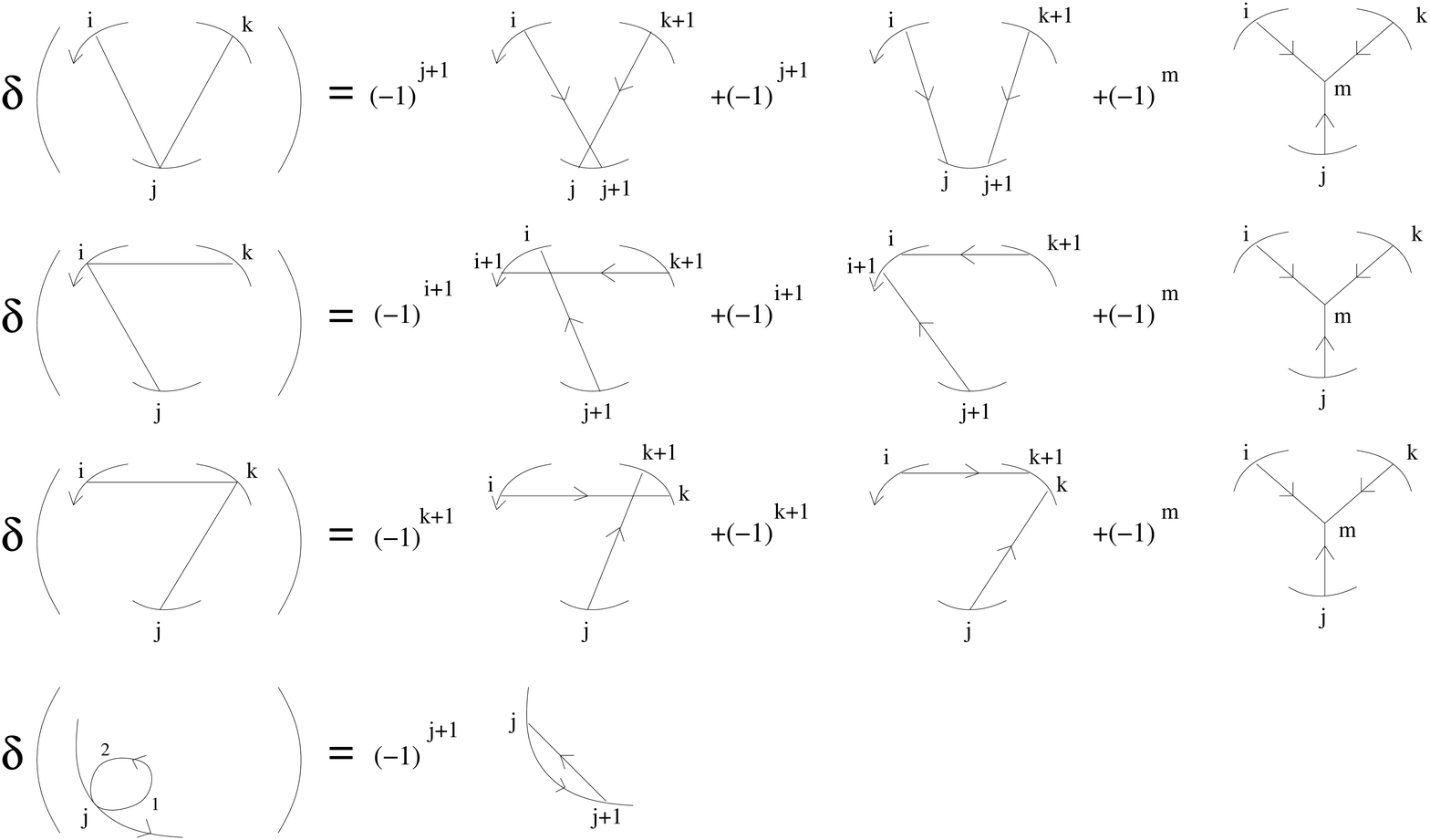}}
\caption{Boundary of a quadrivalent external vertex (odd case)}
\label{de-4ev-d}
\end{figure}

\begin{figure}[h]
\resizebox{5in}{!}{\includegraphics{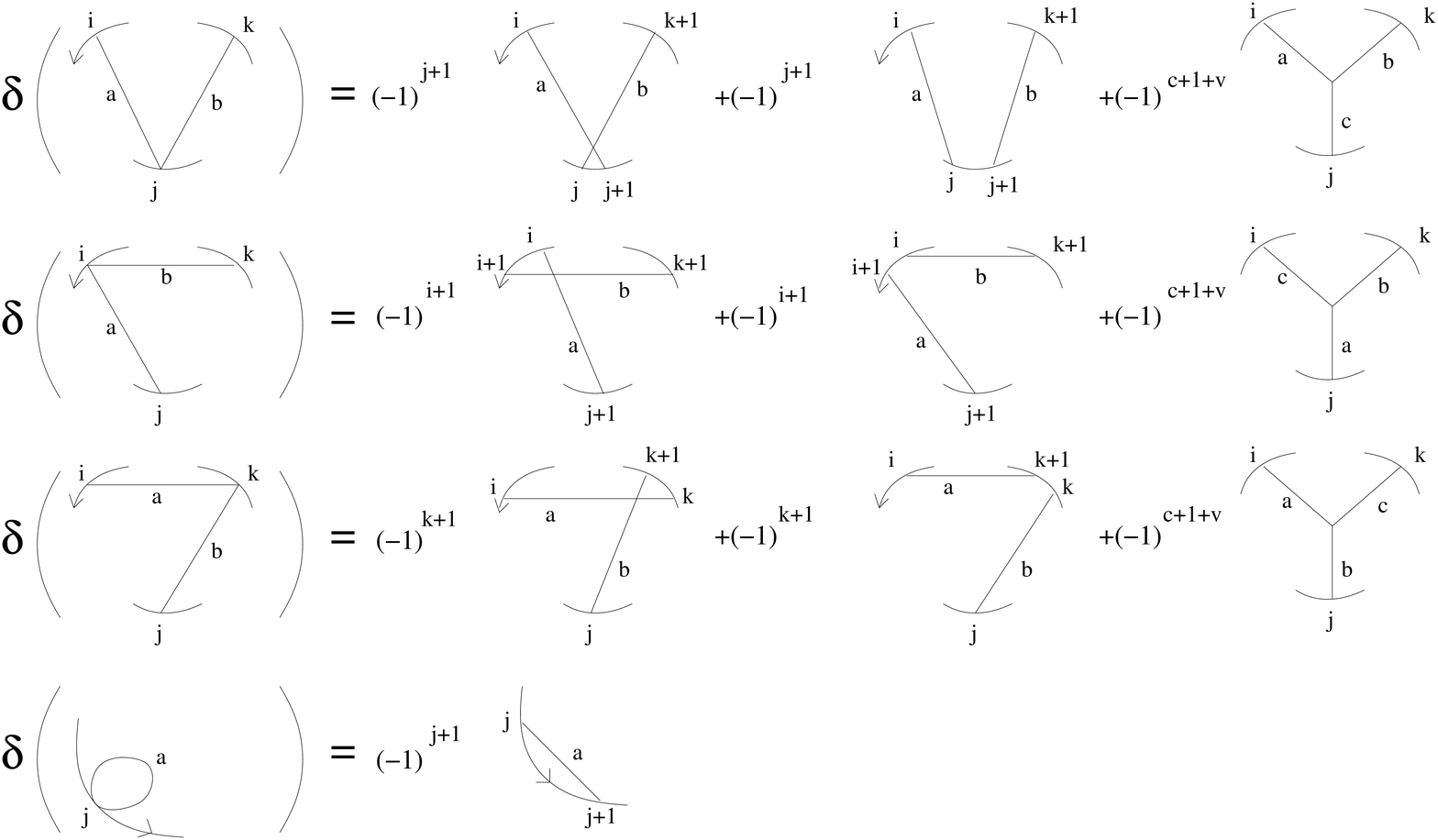}}
\caption{Boundary of a quadrivalent external vertex (even case)}
\label{de-4ev-p}
\end{figure}

For instance, the boundary of a trivalent diagram is always zero, while
the boundary of a diagram $\Gamma$ with all trivalent vertices
except a quadrivalent external vertex is the sum of three
trivalent diagrams as in figure~\ref{de-4ev-d} (odd case) and
\ref{de-4ev-p} (even case).

\begin{figure}[h]
\resizebox{5in}{!}{\includegraphics{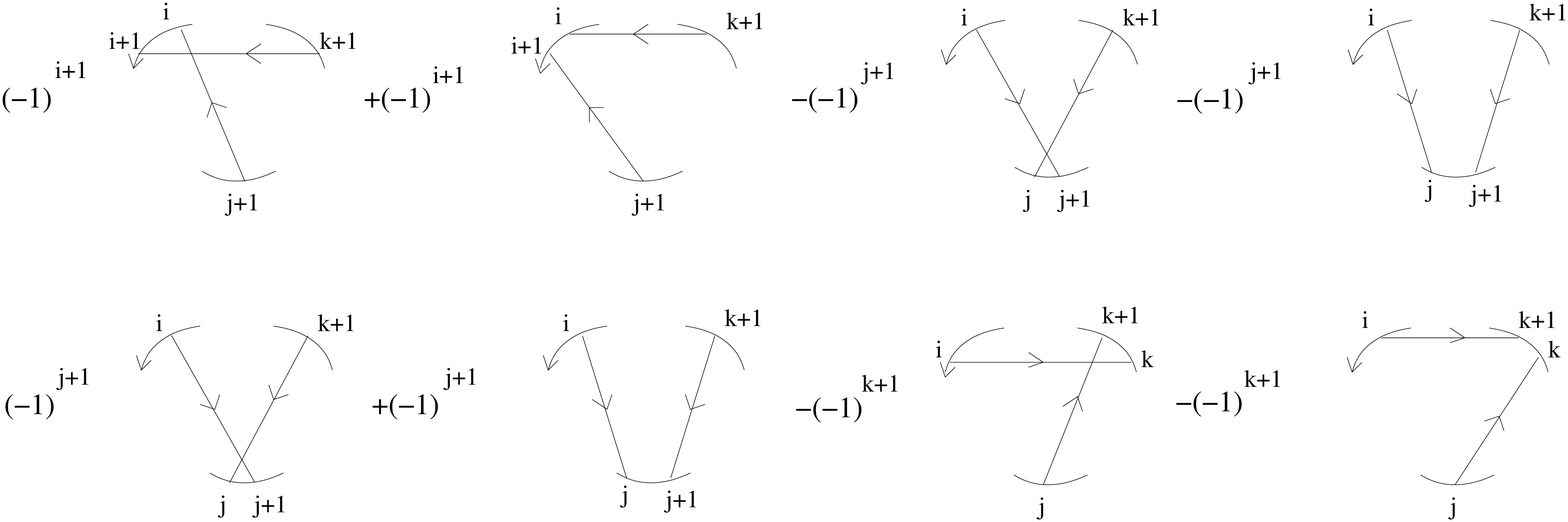}}
\caption{$4T$ relations of odd type}
\label{4t-rel-d}
\end{figure}

\begin{figure}[h]
\resizebox{0.7in}{!}{\includegraphics{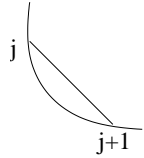}}
\caption{$1T$ relation of odd type}
\label{1t-rel-d}
\end{figure}

\begin{figure}[h]
\resizebox{5in}{!}{\includegraphics{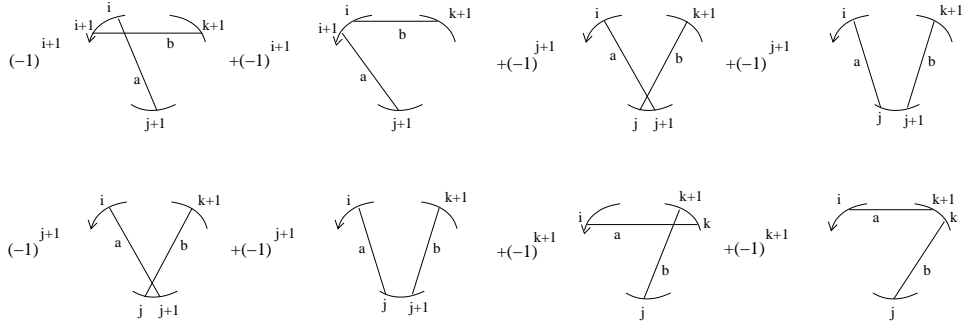}}
\caption{$4T$ relations of even type}
\label{4t-rel-p}
\end{figure}

\begin{figure}[h]
\resizebox{0.7in}{!}{\includegraphics{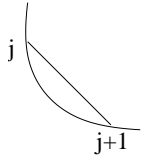}}
\caption{$1T$ relation of even type}
\label{1t-rel-p}
\end{figure}

We claim that all classes of $H_{k,0}(\mathcal D)$ can be
represented by chord diagrams only. More precisely,
$H_{k,0}(\mathcal D_o)$ is isomorphic to the quotient of the space
of chord diagrams $\mathcal {CD}^{k,0}_o$ by the subspace
generated by the diagrams of figure~\ref{4t-rel-d} and 
figure~\ref{1t-rel-d}. Similarly, $H_{k,0}(\mathcal D_e)$ is 
isomorphic to the quotient of the space
of chord diagrams $\mathcal {CD}^{k,0}_e$ by the subspace generated
by the diagrams of figure~\ref{4t-rel-p} and 
figure~\ref{1t-rel-p}. The proof of this fact in
the odd case is given in \cite{BN95}, Thm. 6. The even case is completely
analogous. The subspaces of figure~\ref{4t-rel-d} and \ref{4t-rel-p}
are know as the $4T$ relations, while the subspaces of 
figure~\ref{1t-rel-d} and \ref{1t-rel-p} are the $1T$ relations.

Notice that the $4T$ relations are the boundary of the sum (with
proper sign) of two of the three nontrivalent diagrams considered
in figure~\ref{de-4ev-d} and \ref{de-4ev-p}.

\newpage

\newpage

\section{Cycles and cocycles of imbeddings}
\label{sec:csi}

We now want to recall the definition and main properties of the
{\em configuration space integrals} \cite{BT,CCRL}. 

First, for any compact smooth manifold $M$, we define the open
configuration space $C_q^0(M)$ of $q$ points in $M$ to be the space of
the $q$-uples $(x_1,\ldots,x_q) \in M^q$ such that $x_i\neq x_j$ for 
$i\neq j.$ Then, we consider the Axerlod--Singer--Fulton--MacPherson
compactification $C_q(M)$ of $C_q^0(M)$, given by blowing-up all the
diagonals \cite{Si-c}. These $C_q(M)$ are smooth compact manifolds with
corners. The space $C_q(\R^n)$ is defined as the submanifold of
$C_{q+1}(S^{n})$ where the $(q+1)$st point is fixed to be the point
$\infty\in S^{n}$.

If $\gamma$ is an imbedding of $S^1$ into $\R^n$, then one can force
some of the points in the configuration space to lay on $\gamma$, 
thus obtaining the configuration space $C_{r,s}^0(\R^n,\gamma)$ of $r$
points on $\gamma\subset\R^n$ and $s$ points in $\R^n$. We denote
by $C_{r,s}(\R^n,\gamma)$ its compactification. Putting all these
spaces together, one gets a fiber bundle $p\colon C_{r,s}(\R^n) \to \imbr
n$, whose fibers are compact. 

Next, one considers the maps 
$$
\begin{array}{cccc}
\phi_{ij}\colon & C_{r,s}^0(\R^n) & \to & S^{n-1}\\
 & (x_1,\ldots, x_{r+s}) & \mapsto & \frac{x_i-x_j}{|x_i-x_j|} 
\end{array}
$$
and observes that they extend smoothly to maps $\phi_{ij}
\colon C_{r,s}(\R^n) \to S^{n-1}$. Let $\omega^{n-1}$ be a symmetric top form on 
$S^{n-1}$ which we assume to be concentrated around some fixed direction.
Here ``symmetric'' means that $\omega^{n-1}$ satisfies $\alpha^*\omega^{n-1} = 
(-1)^{n}\,\omega^{n-1}$, where $\alpha$ is the antipodal map on $S^{n-1}$.
The pull-back via $\phi_{ij}$ of $\omega^{n-1}$ is a smooth 
$(n-1)$-form $\theta_{ij}$ on $C_{r,s}(\R^n)$ called 
{\em tautological form}. Forms $\theta_{ii}$ are defined as the
pull-back of $\omega^{n-1}$ via the map
$$
C_{r,s}\stackrel{\pi}\lora C_{r,0}=C_{r}\times \imbr n 
\stackrel{pr_i\times id}\lora S^1\times \imbr n \stackrel{D}{\lora}
S^{n-1}
$$
where $C_r$ is a component of the compactified configuration space of
$r$ points on $S^1$, $\pi$ forgets the $s$ points outside the knot,
$pr_i$ is the projection on the $i$th point and $D$ is the normalized
derivative $D(t,\psi)=\dot\psi(t)/|\dot\psi(t)|$.

We are now ready to define the maps $I^{k,m}$. We consider 
a diagram $\Gamma\in\mathcal D^{k,m}$ and number all its edges and 
vertices (if they are not already numbered by their decoration).
Then to the edge between the vertices $i$ and $j$ we associate the 
tautological form $\theta_{ij}$, and take the wedge product over 
all the edges of $\Gamma$. Finally we integrate (push-forward) this 
product along the map $p\colon C_{r,s}(\R^n) \to \imbr n$, where 
$r$ is the number of internal vertices of $\Gamma$ and $s$ the number 
of external vertices. Extending this map by linearity we obtain
\begin{equation}
\label{chmap}
I^{k,m}\colon \mathcal D^{k,m} \to 
\Omega^{(n-3)k+m}(\imbr n).
\end{equation}
As shown in \cite{CCRL} these are cochain maps for every $n>3$ and 
they induce injective maps in cohomology for $m=0$ .\\ 

One can wonder whether there exists an analogous construction for the
dual theory, namely, whether one can define chain maps
\begin{equation}
\label{chainmap}
(\mathcal D_{k,m}, \de) \to (C_{(n-3)k+m}(\imbr n), \de)
\end{equation}
or, at least, maps in homology
\begin{equation}
\label{homologymap}
H_{k,m}(\mathcal D) \to H_{(n-3)k+m}(\imbr n).
\end{equation}

An answer can be given for the case $m=0$ by associating to
every trivalent chord diagram $\Gamma\in\mathcal {CD}^{k,0}$ a cycle
of $\imbr n$, denoted by $i_k(\Gamma)$, constructed as
follows. Consider first an ``imbedding'' $\psi_\Gamma$ of the diagram
$\Gamma$ in $\R^n$, i.e., an immersion of the oriented circle of
$\Gamma$ into $\R^n$ whose only singularities are transversal double
points, and such that $\psi_\Gamma(t_i) = \psi_\Gamma(t_j)$ if and only
if $t_i$ and $t_j$ are the end-point of a chord in $\Gamma$. We assume
that the indices $i$ and $j$ are the same indices of the chord in $\Gamma$
which has been contracted to form the singular point $\psi_\Gamma(t_i)
= \psi_\Gamma(t_j)$. Next, for every pair of points $t_i, t_j\in S^1$
(say, $i<j$) such that $\psi_\Gamma(t_i) = \psi_\Gamma(t_j)$ and for
every $\mathbf z \in S^{n-3}$, we consider the following loop in
$\R^n$:  
$$
\alpha^{i,j}(\mathbf{z})(t)=
\begin{cases}
0 & \text{if $t\notin [t_i-\epsilon,t_i+\epsilon]$},\\
\mathbf{z}\,\delta\exp
\left(1/[(t-t_i)^2-\epsilon^2]\right)
&\text{if $t\in [t_i-\epsilon,t_i+\epsilon]$},
\end{cases}
$$
with $\epsilon,\delta >0$.
By adding all the loops $\alpha^{i,j}(\mathbf{z})$ to the immersion 
$\psi_\Gamma$ in correspondence with all the double points, we remove 
(blow-up) all the singularities and obtain a family of imbeddings
which is parameterized by $k$ points on $S^{n-3}$, where $k$ is the
number of chords of $\Gamma$. We look at this family of imbeddings as
a $(n-3)k$-cycle in $\imbr n$. 
(We remark that the main result of \cite{CCRL} is that
the cycle of imbeddings obtained from the chord diagram $\Gamma$ is 
dual to the de~Rham cocycle obtained from a trivalent graph cocycle 
containing $\Gamma$, provided this graph cocycle exists). 
Finally we extend these maps $i_k$ by linearity and we have:

\begin{Prop}
\label{triv-4t}
The maps $$i_k\colon \mathcal {CD}^{k,0} 
\to H_{(n-3)k}(\imbr n)$$ defined above, descend to maps on 
$H_{k,0}(\mathcal D)$.
\end{Prop}
\begin{proof}
Since $H_{k,0}(\mathcal D)$ is isomorphic to $\mathcal {CD}^{k,0}$
modulo the $4T$ and $1T$ relations, we have to check that the linear
combination of the diagrams of figure~\ref{4t-rel-d}, \ref{1t-rel-d}, 
\ref{4t-rel-p} and \ref{1t-rel-p} are sent to trivial cycles.

\begin{figure}[ph]
\resizebox{2in}{!}{\includegraphics{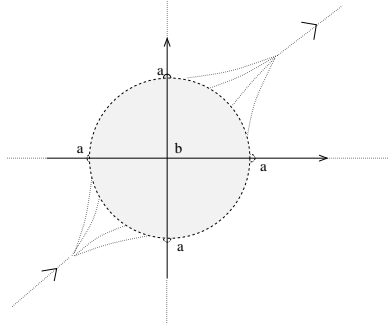}}
\caption{Chain whose boundary gives the $4T$ relation}
\label{ch-4t}
\end{figure}

Consider in fact the four cycles of imbeddings corresponding to 
the diagrams of a $4T$ relation. These cycles (or better, their images)
are equal everywhere except near a small ball in $\R^n$, where they 
are described as follows: two strands of the imbedded circle meet in a
double point $b$ (which is then blown-up) and a third strand is
blown-up around one of the four (half) strands merging in $b$. Denote by 
$i_k(4T)$ the sum of these four cycles. Now, the chain of imbeddings whose 
boundary is $i_k(4T)$ is constructed by considering a chain in $\imbr n$
equal to the cycles in $i_k(4T)$, except that the third strand is lifted 
around $b$ along an $((n-3)+1)$-sphere centered in $b$, with 4 holes 
corresponding to the four (half) strands merging in $b$ 
(see figure~\ref{ch-4t}).

Notice that the orientation of the $((n-3)k+1)$-chain is automatically 
fixed by the orientation of the strand we are lifting and the fact
that we have fixed, once for all, an orientation for $\R^n$.

%
%

As for the $1T$ relations, one can notice that when we collapse and blow-up 
a ``short chord'' (i.e., a chord whose vertices are consecutive vertices 
on the oriented circle) we produce a cycle as in figure~\ref{ch-1t}.

\begin{figure}[ph]
\resizebox{1.7in}{!}{\includegraphics{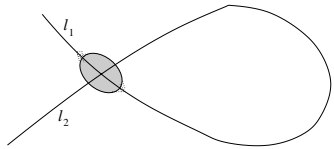}}
\caption{Collapsing and blow-up of a ``short chord''}
\label{ch-1t}
\end{figure}
 
If $n=3$ we have the difference of two ordinary knots which are clearly 
isotopic, namely there exist a 1-chain of imbeddings whose boundary is 
precisely the difference of these two knots. In higher dimensions, one 
can split the ``blow-up'' into two chains, $K_1$ and $K_2$. Suppose in 
fact that the two strands $l_1$ and $l_2$ meeting at the double point 
lie on a plane in $\R^n$; then set $K_1$ to be given by those 
imbeddings lying on one side of the plane determined by $l_1$ and $l_2$, 
and $K_2$ those lying on the other side, in such a way that the cycle 
of imbeddings which ``blows-up'' the double point is $K_1 - K_2$. 
Therefore, the chains $K_1$ and $K_2$ play the role for $n>3$ of the 
two ordinary knots, and it is not difficult to see that there exists 
a higher dimensional analogue to the 1-chain connecting the two 
ordinary knots, namely that there exists a $(n-2)$-chain of embeddings 
whose boundaries are $K_1$ and $-K_2$.
\end{proof}

\section{Proof of Theorem~\ref{thm-main}}
\label{sec:nontriv}

\begin{figure}[ph]
\resizebox{1.7in}{!}{\includegraphics{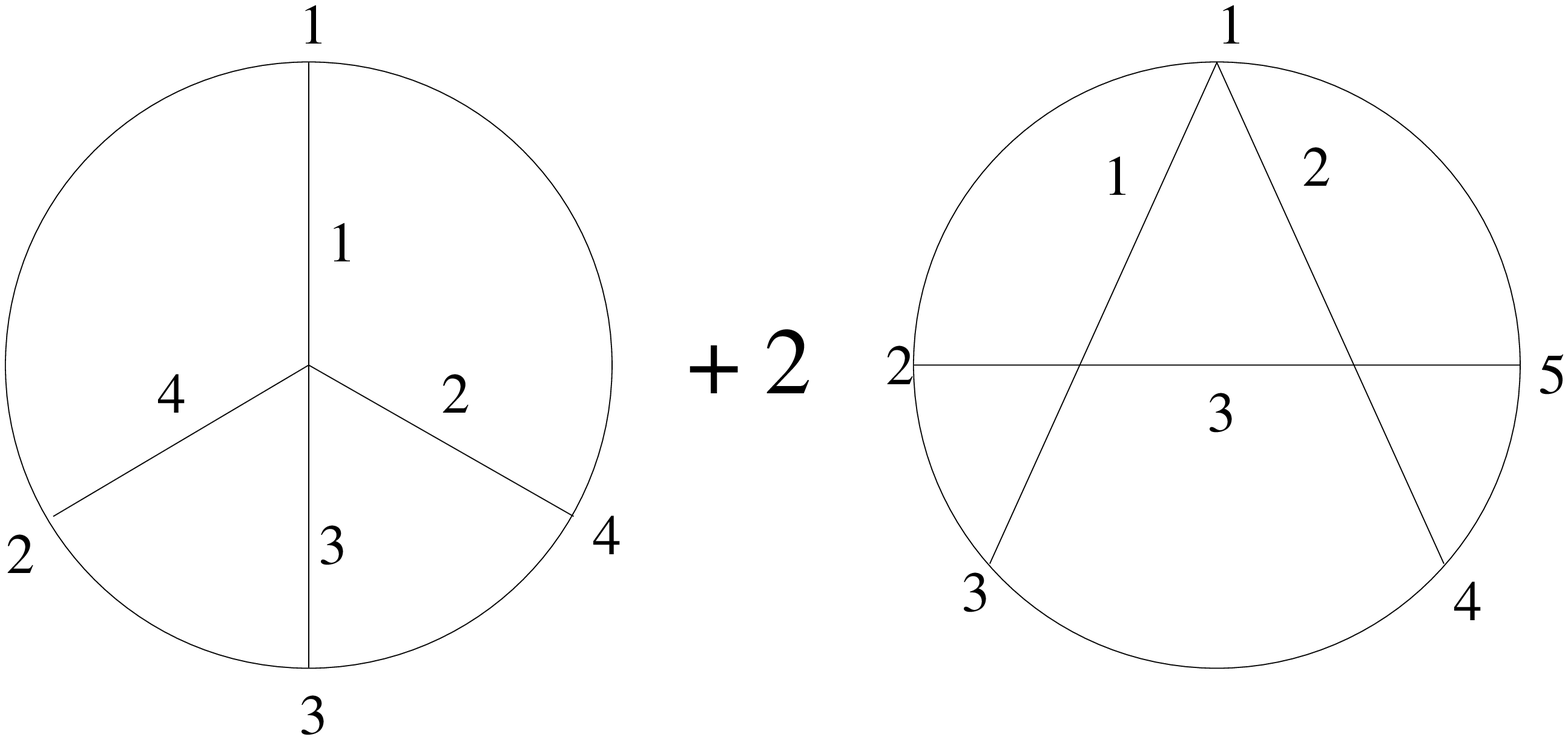}}
\caption{Cocycle of even type of order 3 and degree 1}
\label{pla}
\end{figure}

In the quest for nontrivalent graph cocycles, one easily sees that
there are no nontrivial elements in $H^{1,1}(D_o)$, $H^{2,1}(D_o)$,
$H^{3,1}(D_o)$, $H^{1,1}(D_e)$ and $H^{2,1}(D_e)$, and that the first
example of nontrivial nontrivalent graph cocycle is the generator of
$H^{3,1}(D_e)$. This cocycle is represented in figure~\ref{pla}.

Theorem~\ref{thm-main} is proved by showing that the image of this
graph cocycle through
$$
H^{3,1}(I)\colon H^{3,1}(\mathcal D_e)\to H^{(n-3)3+1} (\imbr{n})
$$
is a nontrivial cohomology class of $\imbr n$. This will follow form
the evaluation of this de~Rham cocycle on the dual cycle, which we
are now going to construct.

\begin{figure}[ph]
\resizebox{5in}{!}{\includegraphics{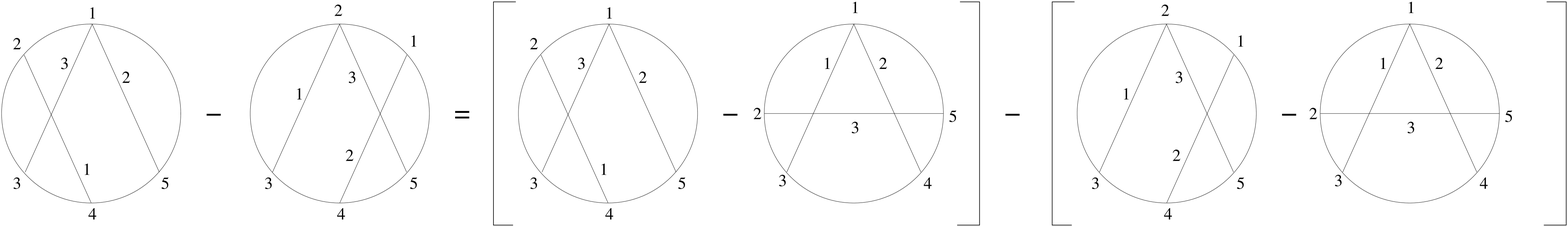}}
\caption{Cycle of even type of order 3 and degree 1}\label{alp}
\end{figure}

Let us first consider the linear combination of diagrams of figure~\ref{alp}. 
We notice that it can be seen as the difference between two chains of 
diagrams $D_1-D_2$ (the $D_i$ are the chains of diagrams in square 
brackets) such that the boundary $\de$ of each $D_i$ gives one of the 
$4T$ relations of figure~\ref{4t-rel-p}. It is also immediate to
check that $D_1-D_2$ is indeed a cycle. Next we associate a 
$((n-3)3+1)$-chain of $\imbr n$ to $D_1$ and $D_2$ as follows: since 
$\de D_i$ is a linear combination of trivalent diagrams, we can 
apply the map $i_3$ of equation~\eqref{homologymap}. As shown in the proof of 
Proposition~\ref{triv-4t} the chain $i_3(D_1)$, resp. $i_3(D_1)$, is the 
boundary of a $((n-3)3+1)$-chain of $\imbr n$ which we denote by the symbol 
$i_3(D_1)$, resp. $i_3(D_2)$. Namely, we have:
\[
i_3(D_i)=\de(i_3(D_i))\qquad i=1,2
\]
More specifically, what we have to do is to consider an ``imbbeding''
of the first diagram of $D_1$, resp. 
$D_2$, i.e., an immersion of the oriented circle of this diagram into 
$\R^n$ which is an imbedding except for a double and a triple point 
obtained by identifying the end-points of each chord. We get $i_3(D_1)$ 
resp. $i_3(D_2)$, by blowning-up these intersections according to the 
constructions of Section~\ref{sec:csi}.

\begin{figure}[ph]
\resizebox{5in}{!}{\includegraphics{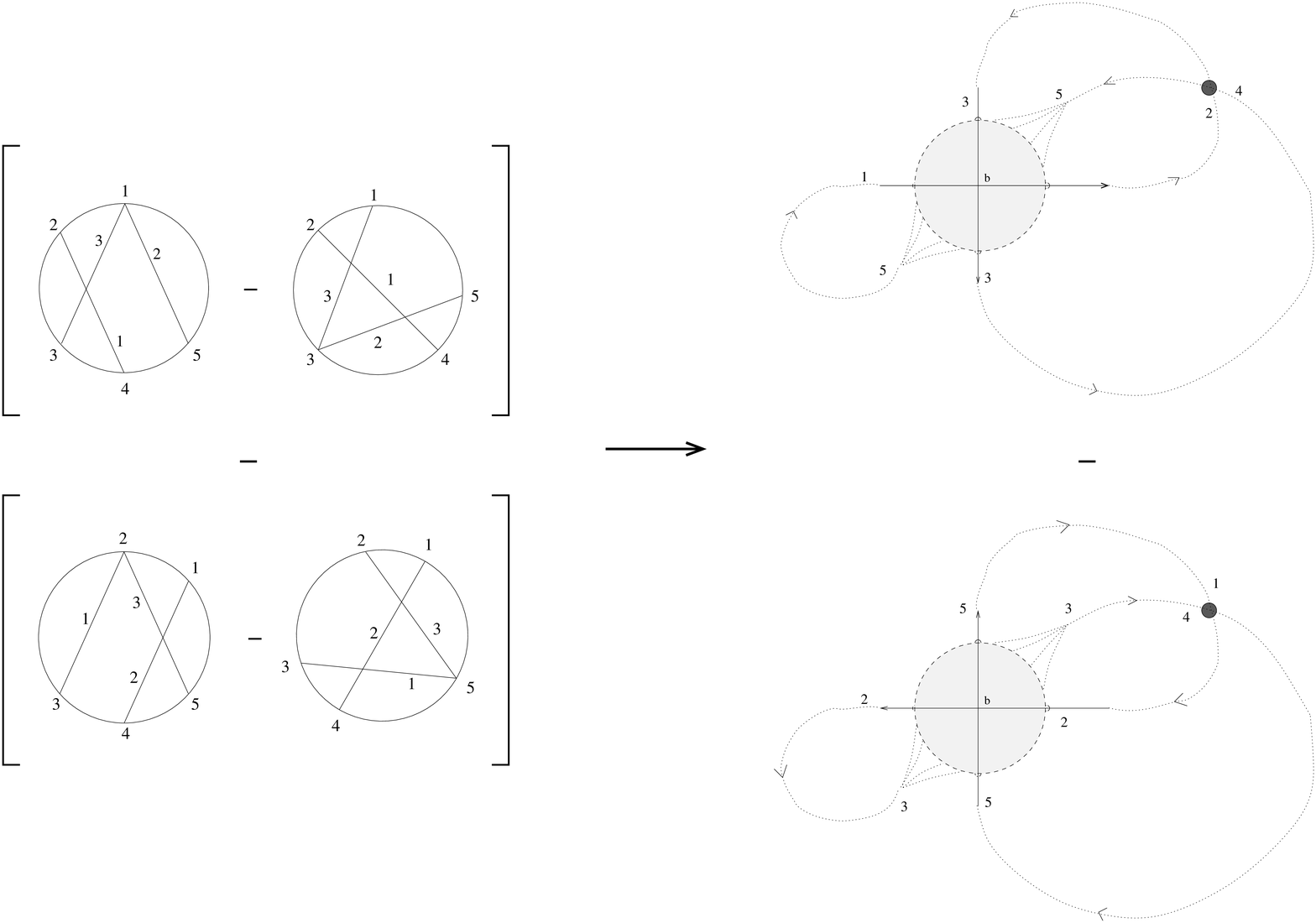}}
\caption{}
\label{ch-alp}
\end{figure}

The difference $i_3(D) = i_3(D_1)-i_3(D_2)$ of these two chains of 
imbeddings is the $((n-3)3+1)$-cycle of $\imbr n$ that we associate to 
$D = D_1-D_2$, and it is represented in figure~\ref{ch-alp}.
In other words, we have extended the map of equation
\eqref{homologymap} to $H_{3,1}(\mathcal D_e)$.

Let us define now
\begin{align*}
\Xi_1=&\int_{C_{4,1}(\R^n)}
\theta_{15}\theta_{45}\theta_{35}\theta_{25}\\
\Xi_2=&\int_{C_{5,0}(\R^n)}\theta_{13}\theta_{14}\theta_{25} 
\end{align*}
so that $\Xi=\Xi_1 + 2\ \Xi_2$ is the 
configuration space integral associated to the diagram of figure~\ref{pla}. 
The Theorem is proved by showing that the evaluation of $\Xi$ on 
$i_3(D_1)-i_3(D_2)$ is different from zero. 

Recall that the tautological forms $\theta_{ij}$ are constructed 
using a symmetric top form $\omega^{n-1}$ on $S^{n-1}$, concentrated around
some fixed direction, say $(1,0,\ldots,0)\in S^{n-1}\subset\R^n$. 
This means that we are requiring $\omega^{n-1}$ to have support near 
$(1,0,\ldots,0)$ and near $(-1,0,\ldots,0)$. 
Moreover we can always suppose that the images in $\R^n$ of the chains 
$i_3(D_1)$ and $i_3(D_2)$ lay on a plane perpendicular to $(1,0,\ldots,0)$, 
except near the blow-ups. 

Therefore the only contributions to the evaluation of $\Xi_1$ and $\Xi_2$ 
on $i_3(D)$, arise from those parts of the configuration space integral 
in which the tautological forms ``points toward the direction 
$(1,0,\ldots,0)$''. This implies in particular that the evaluation of $\Xi_1$ 
on $i_3(D)$ is zero and that, in the evaluation of 
$\Xi_2$ on $i_3(D)$, the only nontrivial contributions
arise from that part of the configuration space integral $\Xi_2$
in which $\theta_{25}$ is integrated near the blow-up of the double
point and $\theta_{14} \theta_{13}$ is integrated near the blow-up of the 
triple point.

Moreover, it not difficult to see that the evaluation of $\Xi_2$ on 
$i_3(D)$ can be explicitely computed as twice the degree of the map
$$
\begin{array}{cccc}
f \colon & C_{5,0}^0(\R^n)    & \lora   & S^{n-1}\times S^{n-1}
\times S^{n-1}\\
& (x_1,x_2,x_3,x_4,x_5)  & \mapsto & \left( \frac{x_1-x_3}{|x_1-x_3|},
\frac{x_5-x_3}{|x_5-x_3|},\frac{x_2-x_4}{|x_2-x_4|} \right).
\end{array}
$$
restricted to a compact subset $K$ of $C_{5,0}^0(\R^n)$. This subset
$K$ is given by the cartesian product of two 1-dimensional manifolds, 
two $(n-2)$-sphere and a $(n-1)$-manifold, imbedded in $\R^n$ as
follows: each of the 1-dimensional manifolds is linked to a
$(n-2)$-sphere, and the $(n-1)$-manifold wraps around one of the
spheres. 

For instance if we consider $x=((1,0,\ldots,0), (1,0,\ldots,0),
(1,0,\ldots,0))$, we can easily see that the point is regular and that
the number of elements in the counter-image is one. In other words, the 
evaluation of $\Xi_2$ on $i_3(D)$ is different from zero, and 
this concludes the proof of Theorem~\ref{thm-main}. 

\section*{Acknowledgements}

This work arises from a conversation with Dylan Thurston. I also thank
Alberto Cattaneo, Paolo Cotta-Ramusino and Victor Tourtchine for
interesting discussions and useful comments.


\begin{thebibliography}{99}

\bibitem{BN95} D. Bar-Natan, On the Vassiliev knot invariants,
{\em Topology} {\bf 34} (1995) 423--472 

\bibitem{BT} R. Bott and C. Taubes, On the self-linking of knots,
{\em J.\ Math.\ Phys.\ } {\bf 35} (1994) 5247--5287 

\bibitem{Bott} R. Bott, Configuration spaces and imbedding invariants,
{\em Turkish J. Math.\ }{\bf 20} (1996) 1--17 

\bibitem{CCRL} A. S. Cattaneo, P. Cotta-Ramusino and R. Longoni, 
Configuration space integrals and Vassiliev classes in any 
dimension, {\em Algebr.\ Geom.\ Topol.\ } {\bf 2} (2002) 949--1000 

\bibitem{CCRRo}  A. S. Cattaneo, P. Cotta-Ramusino and C. Rossi,
Loop observables for BF theories in any dimension and the
cohomology of knots, {\em Lett.\ Math.\ Phys.\ } {\bf 51} (2000) 301--316 

\bibitem{CCRo}  A. S. Cattaneo and C. Rossi,
Higher-dimensional $BF$ theories in the Batalin-Vilkovisky formalism: 
the BV action and generalized Wilson loops, {\em Comm.\ Math.\ Phys.\ } 
{\bf 221} (2001) 591--657

\bibitem{K} M. Kontsevich, Feynman diagrams and low-dimensional topology,
{\em First European Congress of Mathematics, Paris 1992, Volume II,
Progress in Mathematics,} {\bf 120} (Birkh\"auser, Boston, MA, 1994)

\bibitem{Si} D. Sinha, On the topology of spaces of knots, 
\texttt{math.AT/0202287}

\bibitem{Si-c} D. Sinha, Manifold theoretic compactifications of 
configuration spaces, \texttt{math.GT/0306385}

\bibitem{T1} V. Tourtchine, On the homology of the spaces of long knots, 
in {\em Advances in Topological Quantum Field Theory}
NATO Sciences Series by Cluwer (2002)

\bibitem{V} V. Vassiliev, Cohomology of knot spaces, in
{\em Theory of Singularities and its Appl.\ } (ed.\ V.~I.~Arnold) 
Adv.\ in Sov.\ Math.\ {\bf 1} (AMS, Providence, RI, 1990) 23--69


\end{thebibliography}
\end{document}